\documentclass[a4paper, 12pt]{amsart}

\usepackage{amsfonts}                   

\usepackage{latexsym}                   
\usepackage{amssymb}
\usepackage{mathptmx}  
\usepackage{graphicx}
\usepackage[all]{xy}     
\newtheorem{teorema}{Theorem}

\newtheorem{propos}[teorema]{Proposition}
\newtheorem{corol}[teorema]{Corollary}
\newtheorem{ex}{Example}[section]
\newtheorem{rem}{Remark}[section]

\newtheorem{defin}[teorema]{Definition}

\def\defin{\par\ifdim\lastskip<\smallskipamount\removelastskip
  \smallskip\fi\noindent{\bf\ignorespaces
Definition\unskip:\enspace}\rm \ignorespaces}

\def\bit{\begin{itemize}}
\def\eit{\end{itemize}}
\def\be{\begin{equation}}
\def\ee{\end{equation}}
\def\beq{\begin{eqnarray}}
\def\eeq{\end{eqnarray}}

\def\ba{\begin{array}}
\def\ea{\end{array}}
\def\bt{\begin{teorema}}
\def\et{\end{teorema}}
\def\bp{\begin{propos}}
\def\ep{\end{propos}}
\def\bl{\begin{lemma}}
\def\el{\end{lemma}}
\def\bc{\begin{corol}}
\def\ec{\end{corol}}
\def\br{\begin{rem}\rm}
\def\er{\end{rem}}
\def\bex{\begin{ex}\rm}
\def\eex{\end{ex}}
\def\bd{\begin{defin}}
\def\ed{\end{defin}}
\def\demo{\par\noindent{\bf Proof.\ }}
\def\enddemo{\ $\Box$\par\vskip.6truecm}

     \def\nin{\noindent}

\def\ES{\varnothing}  \def\IN{\infty}  

\def\R{{\mathbb {R}}}   \def\a {\alpha} \def\b {\beta}\def\g{\gamma}
     \def\d {\delta} \def\e{\varepsilon}
\def\C{{\mathbb C}}      
                 \def\o{\omega}\def\p{\partial}
                 \def\s{\sigma}
\def\P{{\mathbb P}}

                                                                                                \def\G{\Gamma}
                                                                                                \def\O{\Omega}

  \def\smi{\smallsetminus} \def\ssmi{\!\smallsetminus\!}

\def\sbs{\!\subset\!}

\def\Sbs{\!\Subset\!}

\def\Til {\widetilde}
\def \Hat {\widehat}

\def\oli{\overline}

\def\benu{\begin{enumerate}} \def\eenu{\end{enumerate}}
\def\beqn{\begin{eqnarray*}}  \def\eeqn{\end{eqnarray*}}

\def\beqn{\begin{eqnarray*}}  \def\eeqn{\end{eqnarray*}}

\date{\today}
\begin{document}
\title[Cohomology and removable subsets]{Cohomology and removable subsets}
\author{Alberto Saracco}
\address{A. Saracco: Dipartimento Di Matematica,
Universit\`{a} di Roma \textquotedblleft Tor
Vergata\textquotedblright, Via Della Ricerca Scientifica 1 ---
I-00133 Roma, Italy}
              \email{alberto.saracco@gmail.com}
\author{Giuseppe Tomassini}
\address{G. Tomassini: Scuola Normale Superiore, Piazza dei Cavalieri, 7 --- I-56126 Pisa, Italy}
\email{g.tomassini@sns.it}
\keywords{Cohomology $\cdot$ pseudoconvexity $\cdot$ analytic extension $\cdot$ removable singularities}
\subjclass[2000]{Primary 32D15, 32D20 Secondary 32C25, 32F10}
 \date{\today}
 \begin{abstract} Let $X$ be a (connected and reduced) complex space. A $q$-{\it collar} of $X$ is a bounded domain whose boundary is a union of a strongly $q$-pseudoconvex, a strongly $q$-pseudoncave and two flat (i.e. locally zero sets of pluriharmonic functions) hypersurfaces. Finiteness and vanishing cohomology theorems obtained in \cite{ST1}, \cite{ST2} for semi $q$-coronae are generalized in this context and lead to results on extension problem and removable sets for sections of coherent sheaves and analytic subsets.
 \end{abstract}
 \maketitle
 \tableofcontents
\section{Introduction.}\label{INTR}\nin
Let $X$ be a (connected and reduced) complex space. We
recall that $X$ is said to be {\it strongly} $q$-{\it pseudoconvex} in the sense of
Andreotti-Grauert~\cite{AG} if there exist a compact subset $K\subset X$ and a smooth function
$\varphi:X\to\R$, $\varphi\ge 0$, which is strongly $q$-plurisubharmonic on $X\smi
K$ and such that:
\begin{itemize}
 \item[a)] $0=\min\limits_X\,\varphi<\min\limits_K\,\varphi$;
\item[b)] for every $c>\max\limits_K\,\varphi$ the subset
$$
B_c=\{x\in X:\varphi(x)<c\}
$$
is relatively compact in $X$.
\end{itemize}
If $K=\ES$, $X$ is said to be $q$-{\it complete}. We remark that, for a space, being $1$-complete is equivalent to being Stein.

Replacing the condition b) by
\bit
\item[b')] for every $0<a<\min\limits_K\,\varphi$
and $c>\max\limits_K\,\varphi$ the subset
$$
B_{a,c}=\{x\in X:a<\varphi(x)<c\}
$$
is relatively compact in $X$,
\eit
we obtain the notion of $q$-{\it corona} (see~\cite{AG}, \cite{AT}). A $q$-corona is said to be {\it
complete} whenever $K=\ES$.

The extension problem for analytic objects (basically, sections of coherent sheaves, cohomology classes, analytic subsets) defined on $q$-coronae was studied by many authors (see e.g. \cite{AG}, \cite{FG}, \cite{Se}, \cite{Si}, \cite{SiT}). .

In \cite{ST1}, \cite{ST2} we dealt with the larger class of the {\it semi $q$-coronae} which are defined as follows. Consider a strongly
$q$-pseudoconvex space (or, more generally, a $q$-corona) $X$, and
a smooth function $\varphi:X\to\R$ displaying the
$q$-pseudoconvexity of $X$. Let $B_{a,c}\sbs X$ and let
$h:X\to\R$ be a pluriharmonic function such that
$K\cap\{h=0\}=\ES$. A connected component of
$B_{a,c}\smi\{h=0\}$ is, by definition, a {\it semi}
$q$-{\it corona}. If $X$ is a complex manifold the zero set $\{h=0\}$ can be replaced by a Levi
flat hypersurface. For singular spaces, by definition {\it Levi flatness} means locally zero set of a pluriharmonic function. 

Finiteness and vanishing cohomology theorems proved there lead to results of this type: depending on $q$, analytic objects given near the convex part of the boundary of a {\it semi} $q$-{\it corona} fill in the hole. 

In this paper we consider a more general situation. Let $X$ be a strongly $q$-pseudoconvex space, $C=B_{a,c}=B_c\smi{\overline B}_a$ a $q$-corona. Let $\Sigma_1$, $\Sigma_2$ two Levi-flat hypersurfaces in a neighbourhood of $\oli B_c$ such that 
$$B_c\cap\Sigma_1\cap\Sigma_2=\Sigma_1\cap K=\Sigma_2\cap K=\ES,$$
and $\Sigma_1\cap B_c\neq \ES$, $\Sigma_2\cap B_c\neq\ES$ are nonempty connected subsets. We also assume that $\Sigma_1=\{h_1=0\}$, $\Sigma_2=\{h_1=0\}$ where $h_1$, $h_2$ are pluriharmonic on neighbourhoods $W_1$, $W_2$ of $\Sigma_1\cap B_c$, $\Sigma_2\cap B_c$ respectively.  Let $Q$ be the open subset of $B_c$ bounded by $\Sigma_1\cap B_c$, $\Sigma_2\cap B_c$ and a part of ${\rm b}B_c$. We assume that $Q$ is connected and that $B_c \smi\oli Q$ has two connected components, $B_+$ and $B_-$, and define $C_0=Q\cap C$, $C_+=B_+\cap C$, $C_-=B_-\cap C$. The domain $C_0$ is called a $q$-{\it collar} (see fig.\ 1). A $q$-collar is said to be {\it complete} if $K=\ES$. Note that $C_+$ and $C_-$ are semi $q$-coronae. 
\begin{figure}[h]
	\begin{center}
		\includegraphics[width=0.95\textwidth]{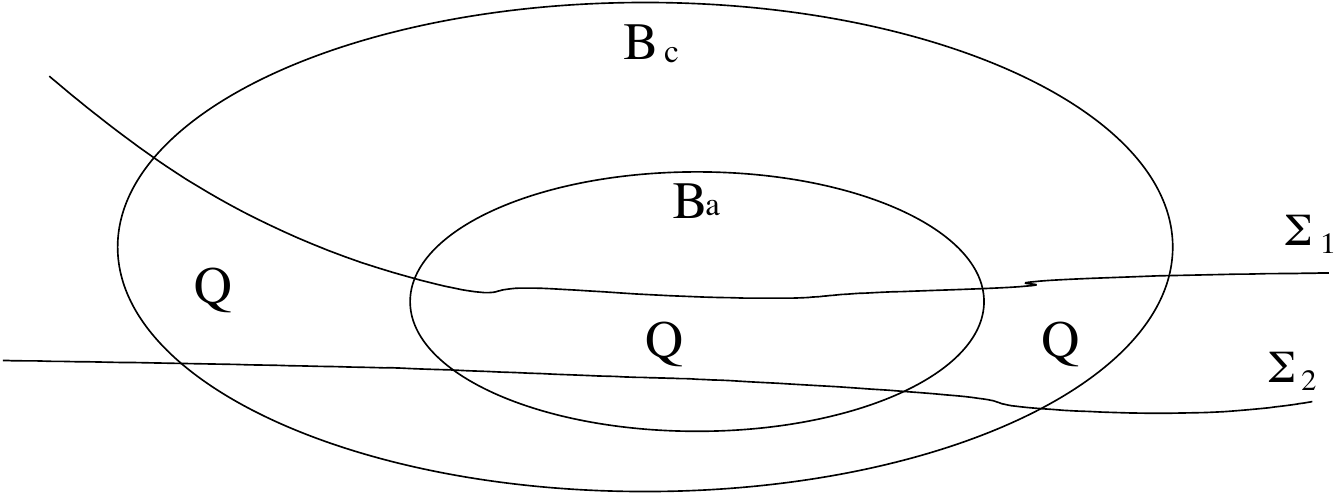}
	\end{center}
	\label{fig:q-collare}
	\caption{A $q$-collar $C_0=Q\cap(B_c\smi \oli B_a)$. In spite of the figure, $C_0$ is connected.}
\end{figure}

Observe that $q$-collar is a difference of two strongly pseudoconvex spaces. Indeed, consider $1/(c-\varphi)$ which is a strongly $q$-plurisubharmonic exhaustion function for $B_c$. We may suppose that the functions $h_1$, $h_2$ are smooth on all of $X$. Moreover, $\psi_1=-\log h_1^2$, $\psi_2=-\log h_2^2$ are plurisubharmonic in $W_1\smi\{h_1=0\}$, $W_2\smi\{h_2=0\}$ respectively. Let $\chi:\R\to\R$ be an increasing convex function such that $\chi\circ\left (1/(c-\varphi)\right )>\psi$ on a neighbourhood of $B_c\smi W_1$. The function 
$$
\Phi_1=\sup\,\left(\chi\circ\left(1/c-\varphi\right),\psi\right)+\frac{1}{c-\varphi}
$$ 
is an exhaustion function for $B_c\smi \{h_1=0\}$ which is strongly $q$-pluri\-sub\-harmonic in $B_c\smi(\{h_1=0\}\cup K)$. In a similar way we construct an exhaustion function $\Phi_2$ for $B_c\smi \{h_2=0\}$ which is strongly $q$-plurisubharmonic in $B_c\smi(\{h_2=0\}\cup K)$. Then the function $\Phi=\sup\,\left(\Phi_1,\Phi_2\right)_{|Q}$ is an exhaustion function for $Q$ which is strongly $q$-plurisubharmonic in $Q\smi K)$. In order to get the conclusion it is sufficient to apply the same argument starting from $B_a$. 

If $\mathcal F\in{\rm Coh}(B_c)$, we define $$p(\mathcal F)=\inf\limits_{x\in B_c}\,{\rm depth}({\mathcal F}_x),$$ the depth of $\mathcal F$ on $B_c$. If $\mathcal F=\mathcal O$, the structure sheaf of $X$, we set $p(B_c)=p(\mathcal O)$.

The results on the cohomology of $q$-collars, generalizing the ones proved in \cite{ST1}, \cite{ST2}, are established in the first part of the paper (see Section \ref{COH}). They are applied in Section \ref{REM} to study removability. Removability for functions was extensively studied by many authors (see e.g. \cite{STOUT}, \cite{LUP94}, \cite{JO}, \cite{CHST}). We are dealing with removability for sections of coherent sheaves and analytic sets. The main results are contained in Theorems \ref{1}, \ref{corea}, \ref{corea2}, \ref{corea3}. 
 \section{Some cohomology}\label{COH}
This section is dealing with cohomology of $q$-collars and some application to extension of sections of coherent sheaves.
\subsection{Closed $q$-collars}
Let $C_0$ be a $q$-collar in a strongly $q$-pseudoconvex space $X$.

\bt\label{Ac}
Let $\mathcal F\in {\rm Coh}(B_c)$. Then, for $q-1\le r\le p(\mathcal F)-q-2$, the homomorphism
$$
H^r(\oli{Q},\mathcal F)\oplus H^r(\oli{C},\mathcal F)\longrightarrow H^r(\oli{C}_0,\mathcal F)
$$
(all closures are taken in $B_c$), defined by $(\xi\oplus\eta)\mapsto\xi_{|\oli{C}_0}-\eta_{|\oli{C}_0}$, has finite codimension.

If $\Sigma_1=\{h_1=0\}$,  $\Sigma_2=\{h_2=0\}$ where $h_1$ and $h_2$ are pluriharmonic functions near $\Sigma_1\cap\oli B_c$ and $\Sigma_2\cap\oli B_c$, respectively, then
$$
\dim_\C\,H^r(\oli C_0,\mathcal F)<\IN
$$
for $q\le r\le p(\mathcal F)-q-2$.
\et
\demo
Consider the Mayer-Vietoris sequence applied to the closed sets $\oli Q$ and $\oli{C}$
\beq\label{suc1}
\cdots &\to& H^r(\oli Q\cup \oli{C},\mathcal F)\to H^r(\oli Q,\mathcal F)\oplus H^r(\oli{C},\mathcal
F)\stackrel{\d}{\to}\\ &\stackrel{\d}{\to}& H^r( \oli C_0,\mathcal F)\to H^{r+1}(\oli Q\cup \oli{C},\mathcal
F)\to\cdots\nonumber
\eeq
$\d(\xi\oplus \eta)=\xi_{|\oli{C}_0}-\eta_{|\oli{C}_0}$. We have
$$
\oli Q\cup \oli C=B_c\smi U,=
$$
where $U=B_a\smi(B_a\cap \oli Q)$. Thus $U$ is $q$-complete and consequently the groups of compact support cohomology $H^{r}_c(U,\mathcal F)$ are zero for $q\leq r\leq p(\mathcal{F})-q$ \cite{AG}.

From the exact sequence of compact support cohomology
\beq
\cdots &\to& H^r_c(U,\mathcal F)\to H^r(B_c,\mathcal F)\to\\
&\to& H^r(B_c\smi U,\mathcal F)\to H^{r+1}_c(U,\mathcal F)\to\cdots\nonumber
\eeq
it follows that
\begin{equation}\label{isomBc-U}
H^r(B_c,\mathcal F)\stackrel{\sim}{\rightarrow} H^r(B_c\smi U,\mathcal F),
\end{equation}
for $q\leq r \leq p(\mathcal{F})-q-1$.

Since $B_c$ is $q$-pseudoconvex,
$$\dim_\C\,H^r(B_c,\mathcal F)<\IN
$$
for $q\le r$ \cite{AG}, and so
$$
\dim_\C\,H^r( B_c\smi U,\mathcal F)<\IN
$$
for $q\le r\le p(\mathcal F)-q-1$.

From (\ref{suc1}) we see that
$$
\dim_\C H^r( B_c\smi U,\mathcal F)=\dim_\C H^r(\oli Q\cup\oli C ,\mathcal F)
$$
is greater than or equal to the codimension of the homomorphism $\delta$. This proves that the image of the homomorphism
$$
H^r(\oli{Q},\mathcal F)\oplus H^r(\oli{C},\mathcal F)\longrightarrow H^r(\oli{C}_0,\mathcal F)
$$
(all closures are taken in $B_c$), defined by $(\xi\oplus\eta)\mapsto\xi_{|\oli{C}_0}-\eta_{|\oli{C}_0}$ has finite codimension provided that $q-1\le
r\le p(\mathcal F)-q-2$, proving the first assertion of the theorem.

If $\Sigma_1=\{h_1=0\}$,  $\Sigma_2=\{h_2=0\}$ are like in the second part of the statement, then, since $K\cap\left(\Sigma_1\cup\Sigma_2\right)=\ES$, $\oli Q$ has a fundamental system of neighborhoods which are $q$-pseudoconvex spaces, thus, by virtue of \cite[Th\'eor\`eme 11]{AG} we have
$$
\dim_\C\,H^r(\oli Q,\mathcal F)<\IN
$$
for $r\ge q$. On the other hand, $\oli C$ is
a $q$-corona, so
$$
\dim_\C\,H^r(\oli C,\mathcal F)<\IN
$$
for $q\le r\le p(\mathcal F)-q-1$ in view of \cite[Theorem 3]{AT}. 

Summarizing, for $q\le r\le p(\mathcal
F)-q-1$ the vector space $H^r(\oli Q,\mathcal F)\oplus H^r(\oli C,\mathcal F)$ has finite dimension and for $q-1\leq r\leq p(\mathcal F)-q-2$ its image in $H^r(\oli C_0,\mathcal F)$ has finite codimension. Thus, for $q\le r\le p(\mathcal F)-q-2$, $H^r(\oli C_0,\mathcal F)$ has finite dimension.
\enddemo
\bt\label{cC}
Assume that $\Sigma_1=\{h_1=0\}$,  $\Sigma_2=\{h_2=0\}$ where $h_1$ and $h_2$ are pluriharmonic functions near $\Sigma_1\cap\oli B_c$ and $\Sigma_2\cap\oli B_c$, respectively, and $\oli Q\cap K=\ES$. Then
$$
H^r(\oli C,\mathcal F)\stackrel{\sim}{\rightarrow} H^r(\oli C_0,\mathcal F)
$$
for $q\le r\le p(\mathcal F)-q-2$ and the homomorphism
\begin{equation}\label{eqA}
H^{q-1}(\oli Q,\mathcal F)\oplus H^{q-1}(\oli C,\mathcal F)\longrightarrow H^{q-1}(\oli C_0,\mathcal F)
\end{equation}
is surjective for $p(\mathcal F)\geq2q+1$.

If $\oli B_+$ is a $1$-complete space and $p(\mathcal F)\ge 3$, the homomorphism
$$
H^0(\oli Q,\mathcal F)\longrightarrow H^0(\oli C_0,\mathcal F)
$$
is surjective.
\et
\demo
By hypothesis $\oli Q$ has a fundamental system of neighborhoods which are $q$-complete spaces, so $H^r(\oli Q,\mathcal F)=\{0\}$ for
$q\le r$ \cite[Th\'eor\`eme 5]{AG}. From (\ref{isomBc-U}) it follows that $H^r(\oli Q\cup \oli C_0,\mathcal F)=\{0\}$ for
$q\le r\le p(\mathcal F)-q-1$. Thus, from the Mayer-Vietoris sequence (\ref{suc1}) we derive the isomorphism
$$
H^r(\oli C,\mathcal F)\stackrel{\sim}{\rightarrow} H^r(\oli C_+,\mathcal F)
$$
for $q \le r \le p(\mathcal F)-q-2$ and that the homomorphism (\ref{eqA}) is surjective if $p(\mathcal F)\geq 2q+1$.

In particular, if $q=1$ and $p(\mathcal F)\ge 3$ the homomorphism
$$
H^0(\oli Q,\mathcal F)\oplus H^0(\oli C,\mathcal F)\longrightarrow H^0(\oli C_0,\mathcal F)
$$
is surjective, i.e.\ every section $\s\in H^0(\oli C_0,\mathcal F)$ is a difference $\s_1-\s_2$ of two sections $\s_1\in H^0(\oli Q,\mathcal F)$, $\s_2\in H^0(\oli C,\mathcal F)$. Since $B_a$ is Stein, the cohomology group with compact supports $H^1_c(B_a,\mathcal F)$
is zero, and so the Mayer-Vietoris compact support cohomology sequence implies that the restriction homomorphism
$$
H^0(\oli B_c,\mathcal F)\longrightarrow H^0(\oli B_c\smi B_a,\mathcal
F)=H^0(\oli C,\mathcal F)
$$
is surjective, hence $\s_2\in H^0(\oli C,\mathcal F)$ is restriction of $\widetilde{\s}_2\in H^0(B_c,\mathcal F)$. So $\s$ is restriction to $\oli C_+$ of $(\s_1-\widetilde{\s}_{2|\oli B_+})\in H^0(\oli Q,\mathcal F)$, and the restriction homomorphism is surjective.
\enddemo
\bc\label{cD}
Let $q=1$ and $p(B_c)\ge 3$. Then every holomorphic function on $\oli C_0$ extends holomorphically on $\oli Q$.
\ec

\subsection{Open $q$-collars}
Keeping the same notations as above consider an open $q$-collar $C_0$. For the sake of simplicity we assume that $B_c$ is $q$-complete. We also assume that $\Sigma_1=\{h_1=0\}$,  $\Sigma_2=\{h_2=0\}$ where $h_1$ and $h_2$ are pluriharmonic functions on open neighbourhoods $U_1$ and $U_2$ of $\Sigma_1\cap\oli B_c$ and $\Sigma_1\cap\oli B_c$, respectively.
\bt\label{oCw}
Let $B_c$ be $1$-complete and $\mathcal F$ a coherent sheaf on $B_c$ with $p(\mathcal F)\ge 3$. Then the homomorphism
$$
H^0(Q,\mathcal F)\longrightarrow H^0(C_0,\mathcal F)
$$
is surjective.
\et
\demo
Let $s\in H^0(C_0,\mathcal F)$. Fix a couple of positive numbers $\varepsilon=(\varepsilon_1,\varepsilon_2)$ small enough such that $\Sigma_{i,\varepsilon_i}$ defined by $\Sigma_{i,\varepsilon_i}=\{h_i=\varepsilon_i\}$ are connected hypersurfaces, $\Sigma_{i,\varepsilon_i}\cap\oli B_c\cap\oli Q\neq\ES$ and $\Sigma_{i,\varepsilon_i}\cap \oli B_c\sbs U_i$, for $i=1,2$.

Consider the open subset $Q_\varepsilon$ of $Q$ bounded by the hypersurfaces $\Sigma_{i,\varepsilon_i}\cap\oli B_c$, and by a part of ${\rm b}B_c$, and set $C_{0,\varepsilon}=Q_\varepsilon\cap C_0$. In view of Theorem \ref{cC} there exists a section $\Tilde s_\varepsilon\in H^0(\oli Q_\varepsilon,\mathcal F)$ which extends $s_{\vert C_{0,\varepsilon}}$. Now observe that the connected component $W$ of $B_c\smi\Sigma_1$ containing $\Sigma_2$ is Stein. So there exists a strongly pseudoconvex domain $\Omega\Sbs W$ such that the domain $D_\varepsilon$ bounded by $\Sigma_{2,\varepsilon_2}\cap\oli B_c$, $\Sigma_2\cap\oli B_c$ and by a part of ${\rm b}B_c$ is relatively compact in $\O$. By Theorem 5 of \cite{ST1} the section $\Tilde s_\varepsilon$ extends on $\O\cap Q$. Thus $s$ extends on $Q_\varepsilon$. In order to conclude the proof we argue as before with respect to the hypersurfaces $\Sigma_{1,\varepsilon_1}$ and $\Sigma_1$.
\enddemo
In particular, we get the extension of holomorphic functions:
\bc\label{oE} If $B_c$ is a $1$-complete space and
$p(B_c)\ge 3$, every holomorphic function on $C_0$ can be
holomorphically extended on $Q$. \ec 

\bc\label{oEE} 
Let $X$ be a
Stein space. Let $\Sigma_1=\left\{h_1=0\right\}\subset X$, and
$\Sigma_2=\left\{h_2=0\right\}\subset X$ be the zero set of two
pluriharmonic functions, and $S$ be a real hypersurface of $X$ with
boundary, such that $S\cap \Sigma_1=b S=b A_1$, $S\cap \Sigma_2=b
S=b A_2$ where $A_1$ is an open set in $\Sigma_1$ and $A_2$ is an
open set in $\Sigma_2$. Let $D\subset X$ be the relatively compact
domain bounded by $S\cup A_1\cup A_2$ and $\mathcal F$ be a
coherent sheaf with ${\rm depth}(\mathcal F)\geq3$. All sections of
$\mathcal F$ on $S$ extend to $D$. 
\ec
\subsection{Finiteness of cohomology}\label{FINCO} 
Results on the cohomology of $q$-collars obtained in the preceding section concern coherent sheaves defined in larger domains. For the applications that we have in mind it is needed to study cohomology of coherent sheaves which are defined just on collars. This can be done by the same methods used in \cite{ST2} for semi $q$-coronae. We briefly sketch the main points of proofs given there focusing on the case $q=1$. The extension for an arbitrary $q$ demands only technical adjustments.
Keeping the same notations as in Section \ref{INTR} let 
$$
C_0=Q\cap (B_c\smi \oli B_a)=Q\cap B_{a,c}=Q\cap\{x\in X:a<\varphi(x)<c\}
$$
be an open $1$-collar of a Stein space $X$ (see fig.\ 1, page \pageref{fig:q-collare}). $Q$ is the subdomain of $B_c$ bounded by the two Levi flat hypersurfaces $\Sigma_1=\{h_1=0\}$, $\Sigma_2=\{h_2=0\}$. $\Sigma_1$ and $\Sigma_2$ are defined on a neighbourhood of ${\oli B}_c$ where $h_1$ and $h_2$ are pluriharmonic functions near $\Sigma_1$ and $\Sigma_2$ respectively. Thus $Q$ is a Stein domain. By $\Sigma^0_1$, $\Sigma_2^0$ we denote the Levi flat parts of ${\rm b}C_0$ and by $F^0_1$, $F^0_2$ the $1$-pseudoconvex and the $1$-pseudoconcave part respectively.  Since $Q$ is Stein, there exist two families of $1$-pseudoconvex hypersurfaces $\left\{\Sigma^\varepsilon_1\right\}$, $\left\{\Sigma^\varepsilon_2\right\}$, $\varepsilon\searrow 0$,  in a neighbourhood of $\oli Q$, with the following properties
\bit
\item[1)] $\Sigma^\varepsilon_1$, $\Sigma^\varepsilon_2$ bound a strip $Q_\varepsilon\subset Q$ and $\Sigma^\varepsilon_1\to \Sigma_1$, $\Sigma^\varepsilon_2\to\Sigma_2$ as $\varepsilon\searrow 0$;
\item[2)] defining $C^\varepsilon_0=Q_\varepsilon\cap B_{a+\varepsilon,c-\varepsilon}$ we obtain an exhaustion $\left\{C^\varepsilon_0\right\}$ of the collar $C_0$.
\eit
Bump lemma and approximation theorem hold for the closed subsets $\oli C^{\,\varepsilon}_0$ with the same proof as in \cite[Lemma 3.3, 3.9]{ST2} and this enables us to the following results. Assume that ${\rm depth}\,\mathcal F_z\ge 3$ for $z$ near to the pseudoconcave part of the boundary of $C_0$; then
\bit
\item[3)] there exists $\varepsilon_0$ sufficiently small such that if $\varepsilon<\varepsilon_0$ the cohomology spaces $H^1(\oli C^+_{\varepsilon},\mathcal F)$ are finite dimensional;
\item[4)] if $\varepsilon<\varepsilon_0$ there exists $\varepsilon_1<\varepsilon$ such that 
$$
H^1( C^+_{\varepsilon'},\mathcal F)\simeq H^1(\oli C^+_{\varepsilon},\mathcal F) 
$$
for every $\varepsilon'\in]\varepsilon_1,\varepsilon[$. 
\eit
3), 4) have an important consequence, namely that for $\mathcal F$ Theorem A of Oka-Cartan-Serre holds in the following form (see \cite[Corollary 4.2]{ST2}:
\bit
\item[5)] if $\e,\e'$ are as in 4), for every compact subset $K$ of $C^+_{\varepsilon'}\smi\{\varphi>c-\e\}$ there exist sections $s_1,\ldots,s_k\in H^0( C^+_{\varepsilon'},\mathcal F)$ which generate $\mathcal F_z$ for every $z\in K$.
\eit
As an application we get the following extension theorem for analytic subsets
\bt\label{OKA2}
Let $X$ be a Stein space, $C_0=Q\cap \left(B_c\smi\oli B_a\right)\sbs X$ be a complete $1$-collar and $Y$ be a closed analytic subset of $C_0$ such that ${\rm depth}(\mathcal O_{Y,z})\ge 3$ for $z$ near $\left\{\varphi=a\right\}$. Then $Y$ extends to a closed analytic subset on $Q$. 
\et
\demo
Taking into account 5) the proof runs as in \cite[Theorem 4.3 and Corollary 4.4]{ST2}.
\enddemo
\section{Removable sets}\label{REM}
The notion of removable sets was originally given with respect to holomorphic function and the removability problem was extensively studied (see e.g. \cite{STOUT}, \cite{LUP94}, \cite{JO}, \cite{CHST}). Here we want to study the same problem with respect to larger classes of analytic objects, namely the classes of sections of coherent sheaves, of cohomology classes and of analytic sets.

Let $X$ be a complex space, $D$ be a bounded domain. Let $\mathcal F$ be a coherent sheaf on a neighbourhood of $\oli D$. A subset $L$ of the boundary ${\rm b}D$ of $D$ is said to be {\it removable} for (the sections of) $\mathcal F$ or for the cohomology classes with value in $\mathcal F$, of a certain degree $r$, if every section 
$s\in \G({\rm b}D\smi L,\mathcal F)$ or cohomology class $\o\in H^r({\rm b}D\smi L,\mathcal F)$ extends by $\tilde s\in \G(\oli D\smi L,\mathcal F)$ or by $\Til\o\in H^r(\oli D\smi L,\mathcal F)$ respectively. 

Similarly, the subset $L$ is said to be {\it removable} for the (respectively, a given) class of analytic subsets if every analytic subset (of a given class of analytic subsets) defined on a neighbourhood of ${\rm b}D\smi L$ extends by an analytic subset of $\oli D\smi L$. 

\subsection{Coherent sheaves}
Given a coherent sheaf $\mathcal F$ on a complex space $X$ let us denote ${Tor}(\mathcal F)$ the torsion of $\mathcal F$; ${Tor}(\mathcal F)$ is the coherent subsheaf of $\mathcal F$ whose stalk at a point $x\in X$ is
$$
{Tor}(\mathcal F)_x=\left\{s_x\in{\mathcal F}_x:\lambda_xs_x=0\ {\rm for\ some} \ \lambda\,\in\mathcal O_x,\,\lambda\neq 0\right\}.
$$
It can be proved (see\ \cite{A}) that the topology of $\mathcal F$ is Hausdorff if and only if $\mathcal F$ has no torsion, i.e.\ ${Tor}(\mathcal F)=\{0\}$. We denote $T(\mathcal F)$ the analytic subset ${\rm supp }\,{Tor}(\mathcal F)$.

Given a bounded domain $D\sbs X$ let $\mathcal A(D)$ be the algebra ${\rm }C^0(\oli D)\cap\mathcal O(D)$ and for every compact $L\sbs\oli D$ let
$$
\Hat L=\left\{z\in\oli D\ :\ |f(z)|\leq\max_L |f|,\forall f\in {\mathcal A}(D)\right\}
$$
be the $ {\mathcal A}(D)$-envelope of $L$. We want to prove the
following 
\bt\label{1} 
Let  $X$ be an $n$-dimensional manifold,
$D$ a bounded pseudoconvex domain in $X$ with a connected smooth
boundary and $L$ a compact subset of $bD$ such that ${\rm b} D\smi
L$ is a connected, nonempty strongly Levi convex hypersurface. Let $\mathcal F$ be a
coherent sheaf on $X$. Assume that: 
\bit
\item[1)] $D\smi \Hat L$ is connected;
\item[2)] ${\rm depth}({\mathcal F_x})\ge 3$
for every $x\in \oli D$; \item[3)] $\dim_\C T(\mathcal F)\cap\oli D\le n-2$
\eit\nin 
Let $U$ be an open neighborhood of $\oli D\smi L$ or $X\smi(D\cup L)$. Then every section of $\mathcal F$ on $U\smi \oli D$ or $U\smi (X\smi D)$ uniquely
extends to a section on $U\smi\Hat L$ or $D\smi\Hat L$. In particular, if $\Hat L=L$ then $L$ is removable for $\mathcal F$.
\et
\demo
The uniqueness is a consequence of the {\it Kontinuit\"atsatz} and of hypothesis 1) and 2). Indeed let $s_1,s_2$ be sections of $\mathcal F$ on $D$ such that $s_1\equiv s_2$ near ${\rm b}D\smi L$. In view of the hypothesis 2), the support of $s_1-s_2$ is an analytic subset $A$ of $D\smi\Hat L$ with no $0$-dimensional irreducible component   
(see \cite[Th\'eor\`eme 3.6 (a), p.\ 46]{BS}). Let $A_1$ be an irreducible component of $A$. Since ${\rm b}D\smi L$ is strongly Levi convex, in view of the {\it Kontinuit\"atsatz} $A_1$ cannot touch ${\rm b}D\smi L$ so $\oli A_1\cap {\rm b}D\equiv\oli A_1\cap ({\rm b}D\smi L)$. Let $x\in A_1$. In view of the hypothesis 1) there exists $f\in{\mathcal A}(D)$ such that $\max_L \vert f\vert<\vert f(x)\vert$. Consider an exhaustion $W_1\Sbs W_2\Sbs\cdots$ by relatively open subsets of $A_1$, $x\in W_1$. By virtue of the maximum principle, for every $k$ there exists a point $x_k\in {\rm b}W_k$ such that $\vert f(x)\vert<\vert f(x_k)\vert$. Then (passing if necessary to a subsequence) we have $x_k\to y\in \Hat L$ as $k\to+\infty$ and consequently $\vert f(x)\vert\le\vert f(y)\vert\le\max_L \vert f\vert$, a contradiction.

We need now to show the existence of the extension. In order to prove the extension we consider just the case that $U$ is an open neighborhood of $\oli D\smi L$ or $X\smi(D\cup L)$ and $\sigma\in \mathcal
F(U\smi \oli D)$, the proof in the other one being similar.  In view of the hypothesis 1), given a point
$x\in D\smi\Hat L$ there exists $f\in{\mathcal A}(D)$, $f=u+iv$, $u$, $v$ real-valued functions,
such that $f(x)=u(x)=1$  $\max_L |f|<1$; in particular $\max_L
|u|<1$. Then, if $\varepsilon>0$ is sufficiently small and
$C=\{u\le1-\varepsilon\}$, we have $C\cap L=\ES$. Let $V$ be an
open neighborhood of $L$ such that $C\cap\oli V=\ES$. Since ${\rm
b}D\ssmi L$ is strongly pseudoconvex, there exists a
pseudoconvex domain $D_1$ with a smooth boundary satisfying the
following properties:
\begin{itemize}
\item[i)] $D\sbs D_1$, ${\oli D}_1\smi D\sbs U$;
\item[ii)] ${\rm b}D_1\cap {\rm b}D\subset V\cap{\rm b}D$;
\item[iii)] ${\rm b}D_1$ is strongly pseudoconvex at the points of ${\rm b}D_1\smi{\rm b}D_1\cap {\rm b}D$.
\end{itemize}
Since $D_1$ is Stein there exists a strongly pseudoconvex $D_2\Sbs D_1$ which contains the compact subset $\oli D\smi V\cap D$ and such that ${\rm b}(D_2\cap D)\smi {\rm b}D\Sbs W$ (see fig.\ 2).
\begin{figure}[h]
	\begin{center}
		\includegraphics[width=0.90\textwidth]{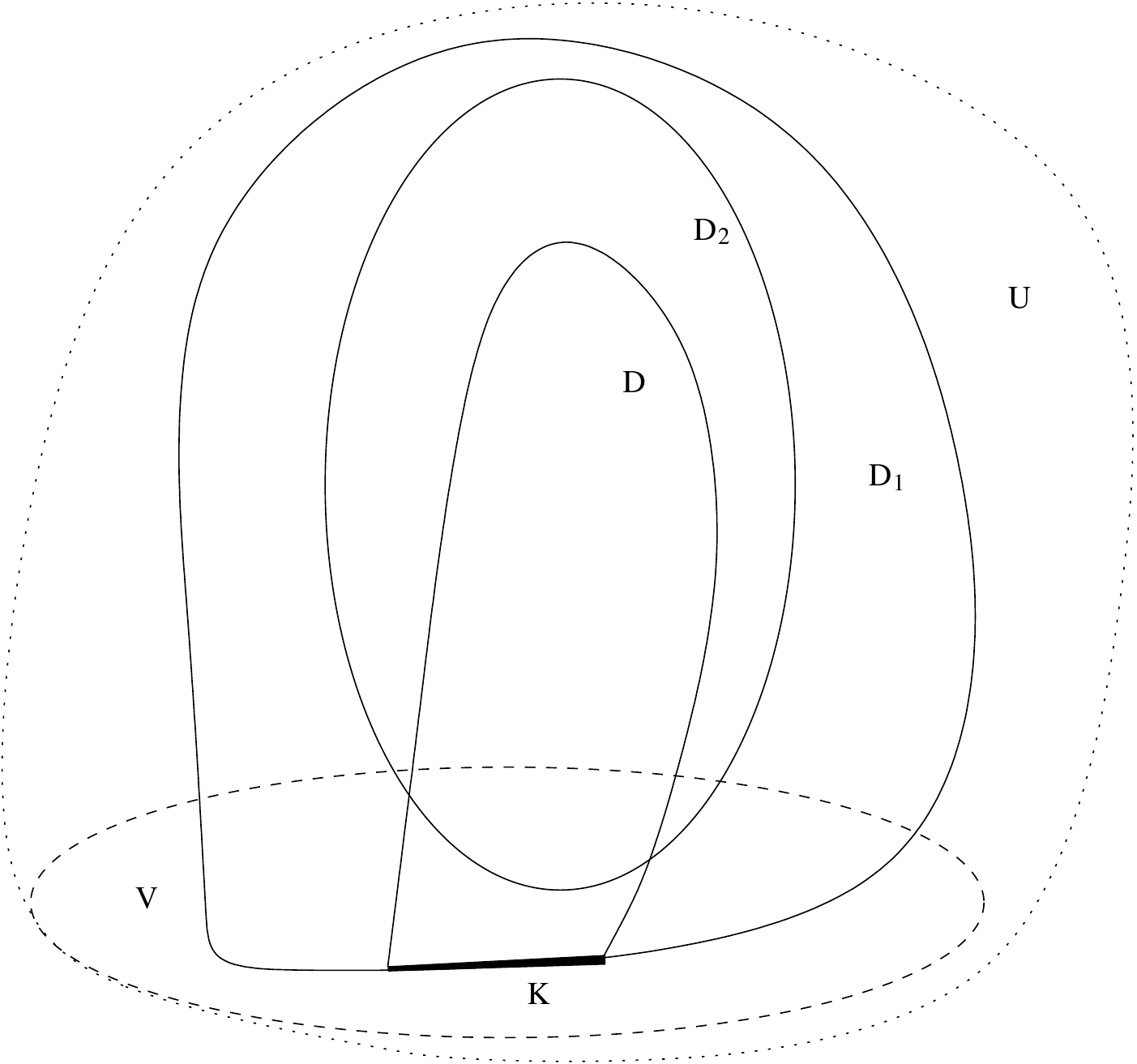}
	\end{center}
	\label{fig:DD1D2}
	\caption{Construction of the three domains $D_1$, $D_2$ and $D_3=D_2\cap D$.}
\end{figure}

The boundary of $D_3=D_2\cap D$ is piecewise smooth but we may regularize it along ${\rm b}D_2\cap {\rm b}D$, thus we may assume that $D_3$ is a smooth strongly pseudoconvex domains $D_3=\{\varrho<0\}$, where $\varrho$ is a strongly plurisubharmonic function on a neighbourhood of $\oli D_3$ and $d\varrho(z)\neq0$ along ${\rm b}D_3$. By the approximation theorem of Kerzman (see \cite{K}) there exists an open neighbourhood $W$ of $\oli D_3$ such that $\mathcal O(W)$ is a dense subalgebra of $\mathcal A(D_3)$. It follows that we may assume that:
\bit
\item[a)] $\sigma\in\mathcal F({\rm b}D_3\cap \{u>1-\varepsilon\})$, where $u$ is pluriharmonic near $\oli D_3$;
\item[b)] $\{u=1-\varepsilon\}$ is smooth and
intersects ${\rm b}D_3$ transversally;
\item[c)] $\{u>1-\varepsilon\}\cap
D_3$ has a finite number of connected components
$D^{(1)},\ldots,D^{(k)}$, which are Stein domains whose boundaries
consist of a part of ${\rm b}D_3$ and of closed subsets
contained in $\{u=1-\varepsilon\}$.
\eit
Moreover, we may suppose that $D^{(i+1)}$
and $D^{(i)}$, $1\le i\le k-1$, are consecutive (i.e. there is a path
$\g\sbs D_3$ joining two points $y'\in D^{(i)}$, $y''\in
D^{(i+1)}$ which does not meet any other connected component $D^{(j)}$ )
and $x\in D^{(1)}$. We denote by $\Sigma _{i}$, $\Sigma _{i+1}$,
$1\le i\le k$, the flat parts of $D^{(i)}$; in particular $\Sigma
_{k+1}=\ES$.

Let us start by $D^{(k)}$. In view of the extension
theorem proved in \cite{ST1}, there exists a unique section
$\sigma_k\in \mathcal F(D^{(k)})$ which extends $\sigma$. Now
consider a positive $\varepsilon'<\varepsilon$ such that the
hypersurface $\{u=1-\varepsilon'\}$ is smooth, intersects ${\rm
b}D_3$ transversally and $D^{(k-1)}$ contains only one connected
component $\Sigma'_k$ of $\{u=1-\varepsilon'\}$. Since
$\O=D^{(k-1)}\cup D^{(k)}\cup \Sigma_k$ is a Stein domain, there
exists a strongly pseudoconvex domain $\O'\Sbs\O$ with the
following properties:
\bit
\item[d)] $\O'$ contains the closed domain bounded by
$\Sigma_k$, $\Sigma'_{k-1}$, ${\rm b}D$, ${\rm b}\O'$ intersects
$\Sigma'_k$ transversally;
\item[e)] no connected component of
$\{u=1-\varepsilon'\}$, $\Sigma'_{k-1}$ excepted, intersects $\O'$.\eit

Thanks again to the quoted extension theorem applied to the
subdomain $\O'$ of $\O$ bounded by ${\rm b}\O$ and intersecting
$\Sigma_k$, we extend $\sigma _k$ by $\sigma' _k$ to $\O'$.
Arguing as above with respect to  the domain bounded by
$\Sigma'_{k-1}$, $\Sigma_{k-1}$, ${\rm b}D$, we extend $\sigma'_k$
to $D^{(k-1)}$ and so on.

In order to finish the proof we have to show that if $\Til\sigma$, $\Til\sigma'$ are two such extensions, defined on $\Til U$ and $\Til U'$ respectively, then $\Til\tau=\Til\sigma'$ on  $\Til U\cap\Til U'$. This is trivially true if $\mathcal F $ is locally isomorphic to a subsheaf of $\mathcal O^N$, in particular if $\mathcal F$ is locally free. 

In our situation consider the difference $\Til\tau=\Til\sigma-\Til\sigma'$ on $\Til U\cap\Til U'$. Since $\mathcal F$ is Hausdorff on $D\smi T$, $T={T}(\mathcal F)$, ${\rm supp}\,\Til\tau\sbs T$. Let $x\in {\rm supp}\,\Til\tau$.  If $B\sbs\Til U\cap\Til U'$ is a sufficiently small Stein neighbourhood of $x$ we have the exact sequences
$$
\xymatrix{&\mathcal O ^p\ar[r]^{\psi}&\mathcal O^q\ar[r]^{\varphi} &\mathcal F\ar[r]& 0}.
$$
$$
\xymatrix{& H^0(B,\mathcal O^p)\ar[r]^{\psi}& H^0(B,\mathcal O^q)\ar[r]^{\varphi}& H^0(B,\mathcal F)\ar[r]&0}.
$$
where $\psi$, $\varphi$ are defined by matrices $(\psi_{ij})$, $(\varphi_{rs})$ of holomorphic functions on $B$. Then $\Til\tau=\varphi(s)$, $s=(s_1,\ldots,s_q)\in H^0(B,\mathcal O^q)$ and $\varphi(s_y)=0$ for every $y\in B\smi T$; consequently
$$
s_{|B\smi T}\in H^0\left(B\smi T,{\sf Ker\,{\varphi}}\right)=H^0\left(B\smi T,{\sf Im\,{\psi}}\right).
$$

It follows that there exist holomorphic functions $g_1,\ldots,g_p$ on $B\smi T$ such that
$$
{s_1}_{|B\smi T}=\sum_{j=1}^p\psi_{1j}g_j,\ldots\dots,{s_q}_{|B\smi T}=\sum_{j=1}^p\psi_{qj}g_j.
$$
Since $\dim_\C T\le n-2$, the functions $g_1,\ldots,g_p$ can be holomorphically extended through $T$ by $\Til g_1,\ldots,\Til g_p$. This implies that $s\in H^0\left(B,{\sf Im\,{\psi}}\right)$, so $s=\psi(\Til g)$, $g=(g_1,\ldots,g_p)$, and consequently $\Tilde\sigma=(\varphi\circ\psi)(\Til g)=0$.

The proof when $U$ is a neighbourhood of $X\smi(D\cup L)$ is similar starting by a pseudoconvex domain $D_1$ with a smooth boundary satisfying the
following properties:
\begin{itemize}
\item[i)] $D_1\sbs D$, $D\smi{\oli D}_1\sbs U$;
\item[ii)] ${\rm b}D_1\cap {\rm b}D\subset V\cap{\rm b}D$;
\item[iii)] ${\rm b}D_1$ is strongly pseudoconvex at the points of ${\rm b}D_1\smi{\rm b}D_1\cap {\rm b}D$.
\end{itemize} 
\enddemo
\br\label{ALEX}
In view of a theorem by Alexander \cite{AL}, condition 1) of Theorem \ref{1} is certainly satisfied if $\Hat L\cap{\rm b}D=L$. Indeed, the connected components $A_i$ of $D\setminus \hat L$ and $B_i$ of $bD\setminus (\hat L\cap bD)$ are in a $1-1$ correspondence given by
$$A_i\leftrightarrow B_i\ \ \ \Longleftrightarrow \ \ \ bA_i\cap bD=B_i.$$
Since $L=\hat L\cap bD$, and $bD\setminus L$ is connected, also $D\setminus \hat L$ is connected.
\er
Condition 1) of Theorem \ref{1} can be dropped also if $L$ is a Stein compact.
\bt\label{corea} Let $X$ be a locally irreducible Stein space, $D$ be
a bounded domain in $X$ with a connected smooth boundary and $L\sbs {\rm b} D$ be a
Stein compact such that ${\rm b} D\smi L$ is
connected. Let $\mathcal F$ be a coherent sheaf on $X$. Assume
that: \bit  \item[1)] ${\rm depth}({\mathcal F_x})\ge 3$ for every
$x\in X$; \item[2)] $ \dim_\C T(\mathcal F)\le n-2$.
\eit\medskip\nin Let $U$ be an open neighborhood of $\oli D\smi
L$. Then every section of $\mathcal F$ on $U\smi \oli D$ uniquely
extends to a section on $U\smi L$. \et
 \demo 
Let 
$$
p(\mathcal F)=\inf\limits_{x\in X}\,{\rm depth}({\mathcal F}_x)
$$  
and $\{U_\alpha\}$ be a fundamental system of Stein neighbourhoods of $L$. Then 
for the compact support cohomology groups we have
$$H^j_c(U_\alpha,\mathcal F)=0,$$
for $j\leq p({\mathcal F})-1$ and every $\a$. Moreover, if $H^j_L(X,\mathcal F)$ denotes the $j^{\rm th}$ local cohomology group with support in $L$ we have the isomorphism
$$
H^j_L(X,\mathcal F)=\varinjlim_{U_\alpha}
H^j_c(U_\alpha,\mathcal F)
$$
(see\ \cite{BS}) hence
$$
H^j_L(X,\mathcal F)=\{0\}
$$
for $j\leq p({\mathcal F})-1$.\\
From the local cohomology exact sequence
$$
\cdots\rightarrow H^j(X,\mathcal F)\rightarrow
H^j(X\setminus L,\mathcal F)\rightarrow H^{j+1}_L(X,\mathcal
F)\rightarrow\cdots\>\>\>,
$$ 
in view of the fact $X$ is a Stein space, we then obtain
$$
H^j(X\smi L,\mathcal F)=\{0\}
$$
for $1\leq j\leq p({\mathcal F})-2$. In particular, since $p({\mathcal F}))\geq3$, we have 
$$
H^1(X\smi L,\mathcal F)=\{0\}.
$$
Let $s\in H^0(bD\smi L,\mathcal F)$. Applying the Mayer-Vietoris sequence to the following closed partition of $X\smi L$
$$
X\smi L=(\overline
D\smi L)\cup[X\smi(D\cup L)]
$$
we get the exact sequence
$$
H^0(\overline D\smi L,\mathcal F)\oplus H^0(X\smi (D\cup L),\mathcal F) \rightarrow H^0(bD\smi L,\mathcal F)\rightarrow H^1(X\smi L,\mathcal F).
$$
Since $H^1(X\smi L,\mathcal F)=\{0\}$ the first homomorphism is onto, so the section $s$ is a difference
$s=s_1-s_2$ of two sections 
$$
s_1\in H^0(\overline D\smi L,\mathcal F),\>\>\>s_2\in H^0(X\smi (D\cup L),\mathcal F).
$$ 
Hence, in order to end our proof, we have to extend the section $s_2$. Consider an open Stein neighbourhood $U$ of $L $. Since, by hypothesis, $p({\mathcal F})\ge 3$, we have $H^1_c(U,\mathcal F)=\{0\}$ and consequently, again from the cohomology exact sequence
$$
H^0(X,\mathcal F)\rightarrow
H^0(X\setminus U,\mathcal F)\rightarrow H^1_c(U,\mathcal F)\rightarrow\cdots\>\>\>,
$$  
we deduce that the homomorphism
$$
H^0(X,\mathcal F)\rightarrow H^0(X\smi\overline U,\mathcal F)
$$
is onto. In particular, there exists a global section  $\tilde s_2$ which extends ${s_2}_{\vert X\smi U}$. Arguing as in the proof of Theorem \ref{1}, we see that $\tilde s_2$ is actually an extension of $s_2$. Thus, $\tilde s=s_1-\tilde s_2$ is a section of $\mathcal F$ on $U\smi L$ which extends $s$. This concludes the proof.
\enddemo
Theorem \ref{corea} can be slightly improved if $X$ is a manifold. Indeed, in that case, under the same hypothesis for $D$, we are allowed to assume that $\mathcal F$ is defined only in a neighbourhood of $\oli D$.

For the proof we need to recall some classical facts about Function Algebra and envelope of holomorphy (see \cite{B}).

Let $X$ be a complex space and ${\mathcal O}(X)$ be the Fr\'echet algebra of all holomorphic functions in $X$. We denote by $\mathcal S(X)$ the {\it spectrum} of ${\mathcal O}(X)$ i.e. the set of all continuous characters $\chi:{\mathcal O}(X)\to\C$ (or, equivalently, the set of all closed maximal ideals of ${\mathcal O}(X)$) equipped with the weak topology. For every $x\in X$ the point evaluation $f\mapsto\d_x(f)=f(x)$, $f\in {\mathcal O}(X)$ is a continuous character and $x\mapsto\d_x$ is a continuous map $i_X:X\to\mathcal S(X)$. Furthermore, for every $f\in{\mathcal O}(X)$ the function $\Hat f:\mathcal S(X)\to\C$ defined by $\Hat f(\chi)=\chi(f)$ is continuous and the set 
$\Hat{{\mathcal O}(X)}=\left\{\Hat f\right\}_{f\in\mathcal S(X)}$ is a subalgebra of ${\rm C}\left({\mathcal S}(X)\right)$. 

Assume now that $X$ is a Stein space. Then, from Oka-Cartan-Serre theory it follows
\bit
\item[$\a$)] $i_X$ is a homeomorphism $X\stackrel{\sim}{\to}\mathcal S(X)$ and there
exists a (unique) complex structure on $\mathcal S(X)$ such that $i_X$ is a biholomorphism
and the dual map $i_X^\ast:{\mathcal O}(\mathcal S(X))\to{\mathcal O}(X)$ is an
isomorphism; in particular $\Hat{{\mathcal O}(X)}={\mathcal O}({\mathcal S}(X))$;
\item[$\b$)] for any complex space $Y$ the functors 
$$
Y\rightarrow {\rm Mor}(Y,X),\>\>\>Y\rightarrow {\rm Hom_{\,cont}}\left({\mathcal
O}(X),{\mathcal O}(Y)\right)
$$\nin
are isomorphic.
\eit
A complex space $X$ is said to have an {\em envelope of holomorphy} if there exists a Stein space $\Til X$ with an open immersion $j:X\hookrightarrow\widehat X$ such that  $j^\ast:{\mathcal O}(\widehat X)\rightarrow {\mathcal O}(X)$ is an isomorphism of Fr\'echet algebras. From the properties ($\a$), ($\b$) it follows that the pair
$(\widehat X,j)$ is uniquely determined (up to isomorphism) by these conditions. Moreover
\bit    
\item[$\g$)] an envelope of holomorphy of a normal space $X$ is also normal (provided
it exists);
\item[$\delta$)] $X$ has an envelope of holomorphy if and only if $\mathcal S(X)$ has a
Stein space structure such that $\left(\mathcal S(X),i_X\right)$ is an envelope of
holomorphy of $X$.
\eit
Cartan, Thullen, Oka and Bishop (see \cite{GR}) proved that for every Riemann domain over $\C^n$,
$p:\O\to\C^n$,  an envelope of holomorphy $\widehat{\O}$ exists  and it is still a domain over
$\C^n$, $\widehat p:\widehat\O\to\C^n$. Using the language of the Function Theory this result can be
stated as follows (see \cite[Theorem\ 2 and Corollary\ 1]{B}: 
\bit
\item[]$\left(\mathcal S(\O),i_\O\right)$ has a complex
manifold structure such that $i_{\O}$ is a holomorphic open immersion and $p=\pi\circ i_{\O}$, and $\Hat{{\mathcal
O}(X)}$ is the algebra of all holomorphic functions in $\mathcal S(\O)$. The natural map
$\pi:\mathcal S(\O)\to\mathcal S(\C^n)\simeq\C^n$ is defined by $\chi\mapsto \left((\chi(z_1),\ldots,\chi(z_n)\right)$ and is
holomorphic of maximal rank. Moreover, $p(\O)\sbs \pi\left(\mathcal S(\O)\right)$. 
\eit

More generally an envelope of holomorphy exists for domains $\O$ over a Stein manifold X
(see \cite{DG}).   

For domains in a Stein space $X$ the envelope of holomorphy could not exist even if $X$ is normal, with isolated
singularities. The first counterexample is due to Grauert (see \cite{DG}). 
\bt\label{corea2}
Let $D$ be a bounded domain of a Stein manifold $X$ with a connected smooth boundary and $L\sbs {\rm b} D$ be a Stein compact  such that 
that ${\rm b}D\smi L$ is connected. Let ${\widehat p}:{\widehat D}\to X$ be the envelope of holomorphy of $D$ and $\mathcal F$ be a coherent sheaf on a neighbourhood $W$ of $\oli{{\widehat p}(\widehat D)}$ satisfying
\bit  
\item[1)] 
${\rm depth}({\mathcal F_x})\ge 3$ for every $x\in W$;
\item[2)] $ \dim_\C T(\mathcal F)\le n-2$.
\eit\medskip\nin 
Let $U\sbs W$ be an open neighborhood of $\oli D\smi L$. Then every section of $\mathcal F$ on $U\smi \oli D$ uniquely extends to $D\smi L$. 
\et
\demo
Let $\widehat W$ be the envelope of holomorphy of $W$, $\widehat p:\Hat W\to X$ be the canonical projection and$j:W\to\Hat W$ be the canonical open embedding of $W$ into $\Hat W$. $j^*:\mathcal O(\Hat W)\to\mathcal O(W)$ is an isomorphism. In particular $\widehat p^*\mathcal F$ is a coherent sheaf on $\widehat W$ with the same depth as $\mathcal F$, which extends $j_*{\mathcal F}$. At this point we argue as in the proof of Theorem \ref{corea}.
\enddemo
\subsection{Analytic sets}\label{ANA}
As for analytic sets, results of removability are obtained arguing as in the proof of Theorem \ref{1} taking into account Theorem \ref{OKA2}. Precisely
 \bt \label{corea3}
Let $X$ be an $n$-dimensional manifold, $D$ be a bounded pseudoconvex domain in $X$ with a connected smooth
boundary and $L$ be a compact subset of $bD$. Assume that: 
\bit
\item[1)] ${\rm b} D\smi L$ is a connected, non-empty strongly Levi convex hypersurface;
\item [2)]  $D\smi \Hat L$ is connected.
\eit
Let $U$ be an open neighborhood of $\oli D\smi L$ and $Y$ be a closed, analytic subset of $U\smi\oli D$ such that ${\rm depth}({\mathcal O_{Y,x}})\ge 3$
for every $x\in U\setminus \oli D$. Then $Y$ extends to an analytic subset $\Til Y$ of $(D\smi\Hat L)\cup U$.
\et
\section{Obstructions to extension} 

The extension theorems proved in the above sections state that, under appropriate conditions, analytic objects like  $CR$-functions, section of coherent sheaves, analytic subsets defined on $bD\smi L$ ($bD\smi L$ being connected) extend---uniquely--- to $D\smi \Hat L$, where $\Hat L$ is the envelope of $L$ with respect to the algebra $\mathcal A(D)$ of holomorphic functions continuous up to the boundary. Natural problems arise about minimality.

In order to state the problem in all generality, given a compct subset $L$ of ${\rm b}D$ we fix a class $\sf C$ of analytic objects and we consider the family $\sf L_{\sf C}$ of all compact subsets $\widetilde L$ of $\overline D$, partially ordered by inclusion, satisfying the following properties
\bit
\item[i)] $\widetilde L\cap bD=L$;
\item[ii)] every analytic object of $\sf C$ defined on $bD\smi L$ extends ---uniquely--- to $D\smi \widetilde L$.
\eit
Suppose that $\sf L_{\sf C}\neq\ES$; then exists in $\sf L_{\sf C}$ some minimal element $L^0_{\sf C}$. One natural problem arises: is $L^0_{\sf C}$ unique? In general, due to polidromy phenomena, the answer could be negative. A second observation is that, at least in the cases already considered, if we have unicity then for the minimal compact $L^0_{\sf C}$ we have the inclusions
$$
L\subset L^0_{\sf C}\subset\Hat L.
$$
The two extremal cases may actually occur. Moreover $L^0_{\sf C}$ heavly depends upon the class $\sf C$. Here are some trivial examples.\medskip
\bit
\item[1)] Let $D=\mathbb B^n\subset\mathbb C^n$ is the unit ball, $L=b\mathbb B^n\cap\{{\mathsf Re}\, z_n\leq 0\}$, $n\geq3$, and $\sf C$ be the class of holomorphic functions.
The minimal compact $L^0_{\sf C}$ is $\Hat L\supsetneq L$.\medskip
\item[2)] Let $D=\mathbb B^n\subset\mathbb C^n$, $L=b\mathbb B^n\cap\{z_2=\cdots=z_n= 0\}=\mathbb S^1\times\{0\}^{n-1}$, $n\geq3$, and $\sf C$ be the class of holomorphic functions.
The minimal compact $L^0_{\sf C}$ is $L \subsetneq\Hat L$.\medskip
\item[3)] Let $D=\mathbb B^n\subset\mathbb C^n$, $L=b\mathbb B^n\cap\{z_{n-2}=\cdots=z_n= 0\}=\mathbb S^1\times\{0\}^{n-1}$, $n\geq5$, and $\sf C_1$ be the class of holomorphic functions, and $\sf C_2$ be the class of analytic sets of codimension $3$.
Then the minimal compacts are
$$L^0_{\sf C_1}\ =\ L\ \subsetneq\ \Hat L\ =\ L^0_{\sf C_2},$$
 as shown by the fact that the analytic set
$$\bigcup_{k\in\mathbb Z} \left\{z_{n-2}=z_{n-1}=0,z_n=\frac1k\right\}$$
does not extend throught $\Hat L$.
\eit
\section{The unbounded case}
Some of the previous results extend to unbounded domains. The following is of particular interest.
\bt
Let $X$ be a complex space and $D$ be a strongly pseudoconvex unbounded domain with a connected boundary. Assume that there exists a sequence $\{p_k\}$ of pluriharmonic functions near $\oli D$ such that 
\bit
\item[1)] $D_k=\left\{x\in D:p_k(x)>0\right\}\subsetneq D_{k+1}=\left\{x\in D:p_{k+1}(x)>0\right\} $;
\item[2)] $D_k\Sbs X$ and $D=\bigcup\limits_{k\ge 1}D_k$.
 \eit
Let $\mathcal F$ be a coherent sheaf on a neighbourhood $U$ of $\oli D$ such that
 \bit 
  \item[3)] ${\rm depth}({\mathcal F_x})\ge 3$ for every
$x\in U$;
 \item[4)] $ \dim_\C T(\mathcal F)\le n-2$.
\eit
Then every section of $\mathcal F$ on $U\smi \oli D$ uniquely
extends to a section on $U$. 
\et
\demo
Fix a section $\s$ of $\mathcal F$ on $U\smi \oli D$. Consider the domain $D_k$. Since $D$ is strongly pseudoconvex, using bump lemma we find a Stein neighbourhood $V_k\sbs U$ of  $\oli D_k$. We may assume that the function $p_k$ is defined on $V_k$, so ${\rm b}D_k\cap {\rm b}D$ is a Stein compact $L_k$ , so we are in position to apply Theorem \ref{corea} and obtain a unique section $\hat\s_k$ of $\mathcal F$ on $V_k\smi L_k$ extending $\s$. Repeating this argument for every $k$, thanks to uniqueness of extension we get the conclusion.
\enddemo
\br
If $X=\C^n$, conditions 1), 2) are implied by the following one
\begin{enumerate}\item[($\star $)]  if $\oli
D^\infty$ denotes the closure of $D\subset \C^n\subset \C\P^n$ in
$\C\P^n$, then there exists an algebraic hypersurface $V$ such
that $V\cap\oli D^\infty=\emptyset$. \end{enumerate}
Under this condition the extension of analytic sets (with discrete singularities) of dimension at least two holds, see \cite{DS}.
\er \medskip

\nin{\bf Acknowledgements}. This research was partially supported by MIUR project \lq\lq Propriet\`{a} geometriche delle variet\`{a} reali e complesse\rq\rq. 

This paper was partly written during a visit of the authors to Seoul National University in November 2007 and was finished thanks to the hospitality given by Scuola Normale Superiore to the first author. Thus we would like to thank SNU and SNS for the hospitality.

\end{document}